\documentclass[11pt]{amsart}

\usepackage{amsmath, amssymb, amscd, mathrsfs, url, pinlabel,setspace,hyperref,verbatim}
\usepackage[margin=1.25in,marginparwidth=1in,centering,letterpaper,dvips]{geometry}
\usepackage{color,dcpic,latexsym,graphicx,epstopdf,comment}
\usepackage[all]{xy}
\usepackage{tikz,pgfplots}
\usepackage{tikz-cd,circuitikz}
\usepackage{enumitem,upquote,tabularx,textcomp,setspace,color}
\usepackage{makecell}
\usepackage{caption} 
\xyoption{rotate}

\newtheorem*{rep@theorem}{\rep@title}
\newcommand{\newreptheorem}[2]{%
\newenvironment{rep#1}[1]{%
 \def\rep@title{#2 \ref{##1}}%
 \begin{rep@theorem}}%
 {\end{rep@theorem}}}
\makeatother

\newtheorem {theorem}{Theorem}
\newreptheorem{theorem}{Theorem}
\newtheorem {lemma}[theorem]{Lemma}
\newtheorem {proposition}[theorem]{Proposition}
\newtheorem {corollary}[theorem]{Corollary}

\numberwithin{equation}{section}
\numberwithin{theorem}{section}

\theoremstyle{definition}
\newtheorem{definition}[theorem]{Definition}
\newtheorem{construction}[theorem]{Construction}
\newtheorem{data}[theorem]{Data}
\newtheorem{notation}[theorem]{Notation}

\newtheorem{remark}[theorem]{Remark}

\newtheorem{example}[theorem]{Example}
\newtheorem*{ack}{Acknowledgement}
\newtheorem*{org}{Organization}

\newlist{pcases}{enumerate}{1}
\setlist[pcases]{
  label=\bf{Case~\arabic*:}\protect\thiscase.~,
  ref=\arabic*,
  align=left,
  labelsep=0pt,
  leftmargin=0pt,
  labelwidth=0pt,
  parsep=0pt
}
\newcommand{\case}[1][]{%
  \if\relax\detokenize{#1}\relax
    \def\thiscase{}%
  \else
    \def\thiscase{~#1}%
  \fi
  \item
}

\newcommand{\bbslash}{\backslash\backslash}
\newcommand{\Z}{\mathbb{Z}}

\newcommand{\rk}{\operatorname{rank}}

\newcommand{\C}{\mathbb{C}}

\newcommand{\bfP}{\mathbf{P}}

\newcommand{\F}{\mathbb{F}}
\newcommand{\Q}{\mathbb{Q}}

\newcommand{\cR}{\mathcal{R}}

\newcommand{\bI}{\mathbb{I}}

\newcommand{\bK}{\mathbb{K}}

\newcommand{\wti}[1]{{\widetilde{#1}}}
\newcommand{\what}[1]{{\widehat{#1}}}

\DeclareFontFamily{U}{mathx}{\hyphenchar\font45}
\DeclareFontShape{U}{mathx}{m}{n}{
      <5> <6> <7> <8> <9> <10>
      <10.95> <12> <14.4> <17.28> <20.74> <24.88>
      mathx10
      }{}
\DeclareSymbolFont{mathx}{U}{mathx}{m}{n}
\DeclareFontSubstitution{U}{mathx}{m}{n}
\DeclareMathAccent{\widecheck}{0}{mathx}{"71}

\newcommand{\pt}{\mathrm{pt}}

\newcommand{\bfi}{\mathbf{i}}

\newcommand{\ina}{I^\natural}
\newcommand{\ish}{I^\sharp}

\usetikzlibrary{calc,intersections}
\tikzset{every picture/.style=thick}
\tikzset{link/.style = { white, double = black, line width = 1.75pt, double distance = 1.25pt, looseness=1.75 }}
\tikzset{crossing/.style = {draw, circle, dotted, minimum size=0.5cm, inner sep=0, outer sep=0}}
\pgfplotsset{compat=1.12}

\usepackage{extarrows}
\usepackage{mathabx}
\usepackage{amsbsy}
\usepackage{appendix}
\usepackage{float}
\usepackage{subfigure}
\usepackage[numbers,sort&compress]{natbib}
\usepackage{amsthm}
\usepackage{geometry}
\usepackage{fancyhdr}
\usepackage{overpic}
\usepackage{mathrsfs}
\usepackage{colonequals}
\usepackage{microtype}
\usepackage{diagbox}
\usepackage{tabularx}

\newcommand{\bpf}{\begin{proof}}
\newcommand{\epf}{\end{proof}}
\newcommand{\bthm}{\begin{theorem}}
\newcommand{\ethm}{\end{theorem}}
\newcommand{\bprop}{\begin{proposition}}
\newcommand{\eprop}{\end{proposition}}
\newcommand{\bcor}{\begin{corollary}}
\newcommand{\ecor}{\end{corollary}}
\newcommand{\blem}{\begin{lemma}}
\newcommand{\elem}{\end{lemma}}
\newcommand{\bdefn}{\begin{definition}}
\newcommand{\edefn}{\end{definition}}
\newcommand{\bcons}{\begin{construction}}
\newcommand{\econs}{\end{construction}}
\newcommand{\bdata}{\begin{data}}
\newcommand{\edata}{\end{data}}
\newcommand{\bexmp}{\begin{example}}
\newcommand{\eexmp}{\end{example}}
\newcommand{\brem}{\begin{remark}}
\newcommand{\erem}{\end{remark}}
\newcommand{\bnot}{\begin{notation}}
\newcommand{\enot}{\end{notation}}
\newcommand{\benu}{\begin{enumerate}}
\newcommand{\benum}{\begin{enumerate}[leftmargin=*]}
\newcommand{\eenu}{\end{enumerate}}
\newcommand{\beq}{\begin{equation}}
\newcommand{\eeq}{\end{equation}}

\newcommand{\ga}{\gamma}

\newcommand{\ep}{\epsilon}

\newcommand{\aand}{~\mathrm{and}~}

\newcommand{\nk}[1]{\nu^{\sharp}_{\mathbb{K}}(#1)}

\definecolor{lygreen}{HTML}{016646}
\title{Singular instanton homology of dual knots}

\author{Fan Ye}
\address{Department of Mathematics\\Harvard University}
\email{fanye@math.harvard.edu}

\makeatletter

\begin{document}

\begin{abstract}
We establish a dimension formula for the unreduced singular instanton homology of dual knots $\widetilde{K}_{p/q}\subset S^3_{p/q}(K)$ for a knot $K\subset S^3$:
$$
\dim I^\sharp(S^3_{p/q}(K),\widetilde{K}_{p/q},\omega; \mathbb{K}) = 2q \cdot r_{\mathbb{K}}(K) + 2|p - q \cdot \nu^\sharp_{\mathbb{K}}(K)|~\mathrm{for}~p/q\neq \nu^\sharp_{\mathbb{K}}(K),
$$where $\omega\subset S^3\backslash K$ is any unoriented $1$-submanifold as the bundle set, $r_{\mathbb{K}}(K)$ and $\nu^\sharp_{\mathbb{K}}(K)$ are integers from the dimension formula of $I^\sharp(S^3_{p/q}(K);\mathbb{K})$ for a field $\mathbb{K}$ defined by Li and the author. In particular, when $\mathbb{K}$ is the two-element field $\mathbb{F}_2$, the reduced singular instanton homology satisfies\[\dim I^\natural(S^3_{p/q}(K),\widetilde{K}_{p/q},\omega;\mathbb{F}_2)=\dim I^\sharp(S^3_{p/q}(K);\mathbb{F}_2)~\mathrm{for}~p/q\neq \nu^\sharp_{\mathbb{F}_2}(K).\]As an application, for a determinant-one knot $K\subset S^3$ other than the unknot and the torus knots $T_{2,3},T_{2,5}$ and a rational $p/q\in (0,6)$ with $p$ odd prime power, the surgery manifold $\widehat{Y}_{p/2q}(\widehat{K})$ is not $SU(2)$-abelian for the double branched cover $\widehat{Y}=\Sigma(S^3,K)$ and the preimage $\widehat{K}\subset \widehat{Y}$ of $K$. We also obtain non-abelian results for $SU(2)$ representations of the knot complement that send the curves of some fixed slope in $(0,6)$ to traceless elements.
\end{abstract}
\maketitle
\section{Introduction}
The reduced and unreduced singular instanton homologies $\ina(Y,K)$ and $\ish(Y,K)$ for a closed oriented connected $3$-manifold $Y$ and an unoriented knot $K\subset Y$ are introduced by Kronheimer--Mrowka \cite{kronheimer2011khovanov} via singular instantons on orbifolds. The ``reduced" and ``unreduced" correspond to the notions in Khovanov homology, as there are spectral sequences from Khovanov homology of a knot in $S^3$ to instanton homology of the mirror knot with the same variant.

In this paper, we study the singular instanton homology of the dual knot in a surgery manifold for a knot $K\subset S^3$, i.e.\ the core of the Dehn filling solid torus. For $p/q\in\Q$, let $S_{p/q}^3(K)$ denote the manifold obtained from $S^3$ by Dehn surgery along $K$ with coefficient $p/q$ and let $\wti{K}_{p/q}\subset S^3_{p/q}(K)$ denote the dual knot. Let $\omega\subset S^3\backslash K$ be an unoriented $1$-submanifold so that it also lies in $S_{p/q}^3(K)$. Baldwin--Sivek \cite{baldwin2020concordance,baldwin2022concordanceII} studied the dimension of the framed instanton homology $\ish(S^3_{p/q}(K),\omega;\C)$ over $\C$. Later, Li and the author \cite{LY2025dimension} extend the results to any coefficient field $\bK$. More precisely, we have the following theorem.
\bthm[{\cite[Theorem 1.1 and Remark 1.2]{LY2025dimension}}]\label{thm: dimension formula}
Suppose $K \subset S^3$ is a knot and suppose $\mu$ is the meridian of $K$. Suppose $p$ and $q$ are co-prime integers with $q\ge 1$. Then there exist a concordance invariant $\nk{K}\in\Z$ satisfying $\nk{\widebar{K}}=\nk{K}$ for the mirror knot $\widebar{K}$. Define \begin{equation*}
	    M=\nk{K}\aand R=r_{\bK}(K)=\min \left\{\dim I^\sharp(S^3_{M}(K);\bK),\dim I^\sharp(S^3_{M}(K),\mu;\bK)\right\}.
	\end{equation*}We have $R=|M|+2h$ for some $h\in\Z_+$ depending on $\bK$. Moreover, we have
	\begin{equation}\label{eq: dim formula rational, main}
	    \dim I^\sharp(S^3_{p/q}(K);\bK)=\dim I^\sharp(S^3_{p/q}(K),\mu;\bK)= qR + |p-qM|.
	\end{equation}
    except possibly when $p/q = M$ and $M$ is even. In the exceptional case, we have
	\begin{equation}\label{eq: differ 2}
	    \left\{\dim I^\sharp(S^3_{M}(K);\bK),\dim I^\sharp(S^3_{M}(K),\mu;\bK)\right\}=\{R,R+2\}.
	\end{equation}
\ethm
\brem\label{rem: SGMME}
Note that for a general choice of $\omega$, we know $I^\sharp(S^3_{M}(K),\omega;\bK)$ is isomorphic to either $I^\sharp(S^3_{M}(K);\bK)$ or $I^\sharp(S^3_{M}(K),\mu;\bK)$. In a recent work of Ghosh--Miller-Eismeier \cite{SGMME}, they could prove that when $\bK=\F_2$, the integers $M_2=\nu^\sharp_{\F_2}(K)$ and $R_2=r_{\F_2}(K)$ are both divisible by $4$. Moreover, they could show that\[\dim I^\sharp(S^3_{M_2}(K);\F_2)=R_2+2\aand \dim I^\sharp(S^3_{M_2}(K),\mu;\F_2)=R_2.\]
\erem
The main theorem of this paper is the following dimension formula for dual knots.
\bthm\label{thm: dual knot dim formula}
Suppose $K\subset S^3$ is a knot and suppose $\omega\subset S^3\backslash K$ is an unoriented $1$-submanifold. Suppose $p$ and $q$ are co-prime integers with $q\ge 1$. Let \[M=\nu^\sharp_{\bK}(K)\aand r_2=r_{\bK}(K)\]be the integers from Theorem \ref{thm: dimension formula}. Then the isomorphism class of $\ish(S_{p/q}^3(K),\wti{K}_{p/q},\omega;\bK)$ does not depend on $\omega$ and\[\dim \ish(S_{p/q}^3(K),\wti{K}_{p/q},\omega;\bK)=\begin{cases}
    2qR + 2|p-qM|& \mathrm{if}~p/q\neq M;\\
    2R ~\mathrm{or}~ 2R+2&\mathrm{if}~p/q=M.
\end{cases}\]When $\bK=\F_2$, the case of $2R+2$ will not happen and\[\dim \ina(S_{p/q}^3(K),\wti{K}_{p/q},\omega;\F_2)=qR + |p-qM|.\]
\ethm
\brem
Based on an on-going project of Bhat, Li, and the author, we expect that when $\bK=\C$, the case of $2R+2$ in Theorem \ref{thm: dual knot dim formula} happens for the right-handed trefoil $T_{2,3}$ and the case of $2R$ happens for the torus knot $T_{2,5}$, or possibly other instanton L-space knots of genera at least $2$. On the other hand, the invariants $\nu^\sharp_\C$ and $r_\C$ were calculated for many knots by Baldwin--Sivek \cite{baldwin2020concordance}.
\erem
\brem
In \cite[\S 5.3]{ghosh2024tau} (see also \cite[Lemmas 4.3 and 4.4]{LY20255surgery} and \S \ref{sec:instanton knot homology} for the identification of $\ina$ and $KHI$), Ghosh--Li--Wong proved that\[\dim \ina(S_{n}^3(K),\wti{K}_{n};\C)=R' + |n-2\tau_I(K)|\]for some integer $R'$ (possibly different from $r_{\mathbb{C}}(K)$) and the instanton tau invariant $\tau_I(K)$. However, the author of this paper is not sure if the dimension formula over $\C$ can be extended to rational slopes.
\erem
\brem
From \cite[Theorem 1.2]{LY2020} and \cite[Theorem 3.20]{LY2021large}, for $p/q\in\Q\backslash\{0\}$, there is a dimension inequality\[\dim \ina(S_{p/q}^3(K),\wti{K}_{p/q};\C)\ge\dim \ish(S_{p/q}^3(K);\C)\]that is indeed from a spectral sequence. However, from Theorems \ref{thm: dimension formula} and \ref{thm: dual knot dim formula}, we have \[\dim \ina(S_{p/q}^3(K),\wti{K}_{p/q};\F_2)=\dim \ish(S_{p/q}^3(K);\F_2)~\mathrm{for}~p/q\neq \nu^\sharp_{\F_2}(K).\]Furthermore, Ghosh--Miller-Eismeier's work mentioned in Remark \ref{rem: SGMME} implies that \[\dim \ina(S_{p/q}^3(K),\wti{K}_{p/q};\F_2)=\dim \ish(S_{p/q}^3(K);\F_2)-2,\]which implies that there cannot be an analogous spectral sequence over $\F_2$ if $\nu^\sharp_{\F_2}(K)\neq 0$. As the construction of the spectral sequence is highly depends on the existence of the Alexander grading and the two kinds of bypass maps with different Alexander grading shifts, the above discussion indicates that either the Alexander grading, or the difference between the two bypass maps does not exist over $\F_2$.
\erem
The singular instanton homology of knots is closely related to $SU(2)$ representations of fundamental groups, in particular for the knot complements whose images of the meridian are traceless. To state the results more precisely, we introduce the following definitions.
\bdefn
Suppose $K$ is a framed knot in a closed oriented connected $3$-manifold $Y$. We write $Y\bbslash K=Y\backslash\operatorname{int}N(K)$ for the knot complement and write $(\mu,\lambda)$ for the meridian-longitude basis of $\partial (Y\bbslash K)\cong T^2$ with $\mu\cdot \lambda=-1$. Recall that $Y$ is called \emph{$SU(2)$-abelian} if all representations\[\rho:\pi_1(Y)\to SU(2)\]have abelian images, i.e.\ factor through a copy of $U(1)\subset SU(2)$. 

If all meridian-traceless $SU(2)$ representations\begin{equation}\label{eq: traceless condition}
    \rho:\pi_1(Y\bbslash K)\to SU(2)~\mathrm{with}~\rho(\mu)=\bfi
\end{equation}have abelian images, then $(Y,K)$ is called \emph{meridian-traceless $SU(2)$-abelian}. Here $\bfi\in SU(2)$ denotes the diagonal matrix with entries $i,-i$. We know any traceless element in $SU(2)$ is conjugate to $\bfi$ and any $SU(2)$ representation $\rho'$ with $\rho'(\mu)$ traceless is conjugate to one satisfying \eqref{eq: traceless condition}. 

Similarly, if all representations\begin{equation*}
    \rho:\pi_1(Y\bbslash K)\to SU(2)~\mathrm{with}~\rho(\mu^p\lambda^q)=\bfi
\end{equation*}have abelian images for some co-prime integers $p$ and $q$, then $(Y,K)$ is called \emph{$p/q$-traceless $SU(2)$-abelian}. Note that $\rho(\mu^p\lambda^q)=\bfi$ implies that $\rho(\mu^{-p}\lambda^{-q})=-\bfi$, so $p/q$-traceless condition is equivalent to $(-p)/(-q)$-traceless condition, which justifies the rational notation.
\edefn
Recall from \cite{BLSY21,farber2024fixed}, the only instanton L-space knots of genera at most $2$ are the unknot $U$ and the torus knots $T_{2,3}$ and $T_{2,5}$. We have the following theorems about $SU(2)$-representations for other knots.
\bthm\label{thm: main SU(2) 1}
Suppose $K\subset S^3$ is a knot that is not $U,T_{2,3},T_{2,5}$ and suppose $\det(K)=|\Delta_K(-1)|$ is the determinant of $K$, where $\Delta_K(t)$ is the Alexander polynomial of $K$. Then for any slope $p/q\in (0,6)$ such that $p$ is an odd prime power and $p$ does not divides $\det(K)$, we know that $(S^3,K)$ is not $p/q$-traceless $SU(2)$-abelian.
\ethm
\bthm\label{thm: main SU(2) 2}
Suppose $K\subset S^3$ is a knot that is not $U,T_{2,3},T_{2,5}$. Suppose $\what{Y}=\Sigma(S^3,K)$ is the double branched cover of $K$ and $\what{K}\subset \what{Y}$ is the preimage of $K$. If $\det(K)=1$, or equivalently $\what{Y}$ is an integral homology sphere, then for any slope $p/q\in (0,6)$ such that $p$ is an odd prime power, we know that $\what{Y}_{p/2q}(\what{K})$ is not $SU(2)$-abelian.
\ethm
\brem
The condition for $p$ is from some non-degeneracy condition \cite[\S 4]{BLSY21}. Since $\det(K)$ does not depend on the choice of $p/q$, there is only finitely many choices of $p$ that divide $\det(K)$. For the cases of $T_{2,3}$ and $T_{2,5}$, one may compute the $SU(2)$ representations directly by \cite[\S 5.1]{HHK2014pillowcase} and \cite{sivek2022menagerie}.
\erem
\brem
Together with Ghosh--Miller-Eismeier's results \cite{SGMME} mentioned in Remark \ref{rem: SGMME}, we could replace the assumption $p/q\in (0,6)$ in Theorems \ref{thm: main SU(2) 1} and \ref{thm: main SU(2) 2} by $p/q\in(0,8)$.
\erem
\begin{org}
In \S \ref{sec:instanton knot homology}, we we review many kinds of instanton homologies and study their relation and their dependence of the bundle sets. In \S \ref{sec: Dimension formula}, we prove Theorem \ref{thm: dual knot dim formula} by considering the integral case and the rational case separately. In \S \ref{sec: SU(2)-representations}, we describe the relation between singular instanton homology and $SU(2)$ representations, and then prove Theorems \ref{thm: main SU(2) 1} and \ref{thm: main SU(2) 2}.
\end{org}
\begin{ack}
The author thanks Zhenkun Li for the discussion about the work \cite{LY2025dimension}, which motivates this project. The author thanks Sudipta Ghosh and Mike Miller Eismeier for sharing the draft of their paper \cite{SGMME}. The author is partially supported by Simons Collaboration \#271133 from Peter Kronheimer and is also grateful to Yi Liu for the invitation to BICMR, Peking University while this project is developed. 
\end{ack}

\section{Instanton knot homology}\label{sec:instanton knot homology}
In this section, we review many kinds of instanton homologies for knots and their relation. We also study the dependence of the bundle sets carefully, similar to the case of framed instanton homology as in \cite[\S 2]{LY2025dimension}.

Recall from \cite[\S 2.1]{LY2025torsion} that there are three approaches to define the framed instanton homology $\ish(Y,\omega)$ for a closed connected oriented $3$-manifold $Y$ and an unoriented $1$-submanifold $\omega\subset Y$ (called the \emph{bundle set}). We extend the discussion to an unoriented knot $K$ in $Y$ that is disjoint from $\omega$ and describe the three approaches to define the \emph{(reduced) instanton knot homology} $\ina(Y,K,\omega)$, with the correspondence\begin{equation}\label{eq: Ina Ish}
    \ish(Y,\omega)=\ina(Y,U,\omega)
\end{equation}for the unknot $U\subset Y$ that bounds a trivial disk. We also add the fourth approach which is only available over $\F_2$ or coefficient rings with $4=0$.

\benu
  \item From \cite[\S 4.3]{kronheimer2011khovanov}, we take the singular instanton homology\[I(Y,K\# H,\omega\sqcup {\alpha}),\]where $H$ is a standard Hopf link in $S^3$ and $\alpha$ is a standard arc connecting two components of $H$. Here we need to pick a basepoint on $K$ to specify the position of connected sum and a tangent vector at the basepoint to specify the arc $\alpha$. Note that this is the original definition of $\ina(Y,K,\omega)$ when $\omega=0$ (the emptyset).
  \item From \cite{floer1990knot} and \cite[\S 5.4]{kronheimer2011khovanov}, we take the original instanton homology \[I(T^3_K,\omega\sqcup S^1_3)^{\psi}.\]Here we write \[T^3=S^1_1\times S^2_1\times S^3_1\aand T^3_K=(Y\bbslash K)\cup (T^3\bbslash S^1_1)\]via gluing the longitudes of $K$ (for some chosen framing) to the meridians of $S^1_1$ and vice versa. The superscript $\psi\subset H^1(T^3_K;\F_2)$ is the 2-element subgroup generated by the torus $R=S^1_1\times S^1_2\subset T^3_K$, and we take the quotient of the usual configuration space by this extra $2$-element group when constructing the instanton homology.
  \item From \cite[\S 7.6]{kronheimer2010knots}, we take the sutured instanton knot homology \[KHI(Y,K,\omega)=SHI(Y\bbslash K,\ga_K=\mu\cup(-\mu),\omega),\]for an oriented meridian $\mu$ of $K$. Here we need to choose two basepoints on $K$ to specify the positions of the meridians. This is only defined over $\C$ and on the homology level since the construction uses the generalized eigenspace decomposition of the $\mu$ operators on the ordinary instanton homology. It is only an isomorphism class from \cite[\S 7]{kronheimer2010knots}, or a projectively transitive system by Baldwin--Sivek's naturality result \cite[\S 9]{BS2015naturality}.
  \item From \cite{kronheimer2019web1} and \cite[\S 2.5]{KM2021Barnatan}, we take the singular instanton homology for webs\[I(Y,K\#\Theta,\omega),\]where $\Theta$ is a standard theta graph in $S^3$. This is only defined over $\F_2$ or coefficient rings with $4=0$ because of the existence of trivalent vertices indicates that the square of the differential $\partial^2$ is a multiple of $4$. Here we need to pick a basepoint on $K$ and an oriented basis $(\ep_1,\ep_2,\ep_3)$ of tangent vector at the basepoint with $\ep_1$ pointing along $K$ to specify the position of the connected sum.
\eenu

All definitions can be generalized to a link $L\subset Y$ if we choose one component of $L$ to place the basepoint or carry out the gluing construction (in the second approach, we need to use singular instanton homology and in the third approach, we need to use sutured instanton homology for singular tangles by Xie--Zhang \cite{XZ2025tangle}). If we take $L=K\sqcup U$ for a split unknot $U$ with basepoint, then the corresponding instanton homology\[\ina(Y,K\sqcup U,\omega)=\ish(Y,K,\omega)\]is called the \emph{unreduced instanton knot homology}. Note that $K$ could be empty in the unreduced variant, for which we will omit it from the notation and obtain the framed instanton homology $\ish(Y,\omega)$. When $\omega=0$, we also omit it from the notation.

From \cite{KM2019deformation,KM2021Barnatan}, there are further constructions about local coefficients for the fourth approach. Moreover, from \cite{kronheimer2011khovanov}, in the first and the fourth approaches, we allow $\omega$ to include arcs with endpoints on $K$. Indeed, we may also extend arcs into circles in the second and the third approaches and allow $\omega$ to include arcs. However, we will not use those extended constructions and the corresponding results for the main results of this paper.

The constructions in the first two approaches inherit $\Z/4$ homological grading, while ones in the last two only inherits a $\Z/2$ homological grading. Here the homological grading can be absolute or relative, depending on $\omega$ and the approaches.

From \cite[\S 5.4]{kronheimer2011khovanov}, the first two definitions are isomorphic via Floer's excision cobordism, and hence the isomorphism intertwines with any cobordism map supported in $Y\bbslash K$ naturally. The first and the last definitions are isomorphic via another excision cobordism \cite[\S 3.3]{KM2019deformation}, and hence also natural for cobordism maps supported in $Y\bbslash K$.

Also from \cite[\S 5.4]{kronheimer2011khovanov}, the second definition is isomorphic to the third definition via a \emph{special} choice of closure in the construction of $KHI(Y,K,\omega)$, which is exactly $(T^3_K,R,\omega\sqcup S^1_3)$. Since $g(R)=1$, by \cite[Theorem 2.5]{baldwin2019lspace} (see also \cite[Lemma 4]{Froyshov2002equi}), we know $\mu(\pt)^2=4$ for the $\mu$ map in the construction of $KHI(Y,K,\omega)$, and then the isomorphism can be made to intertwine with the cobordism map. Indeed, as in \cite[\S 8]{LY2025dimension}, one can replace $\C$ by any field $\bK$ with $\operatorname{char}\bK\neq 2$ and the isomorphism still holds.

Note that Baldwin--Sivek's naturality result \cite{BS2015naturality} only works for closures of genera larger than one, while the relation between closures of genus one and larger genera remains to be isomorphism rather than \emph{canonical} isomorphism. Hence, the best result we can state about the cobordism map is that, for any \emph{fixed} closure and any \emph{fixed} isomorphism between the special closure above and the fixed closure, the isomorphism intertwines the cobordism map supported in $Y\bbslash K$. Similar isomorphism results hold for the unreduced variant by Floer's excision theorem (for sutured version, see \cite[Remark 7.8]{XZ2025tangle}). All isomorphisms above respect the homological gradings that defined on both sides.

Then we list some relation between $\ina$ and $\ish$ as follows, though we only use some of them. Let $\cR$ denote a general coefficient ring and let $\bK$ denote a general coefficient field.

From \cite[Lemma 7.7]{kronheimer2019web1}, we have \beq\label{eq: 2dim}
    \dim \ish(Y,K,\omega;\F_2)=2\dim \ina(Y,K,\omega;\F_2).
\eeq
The equation might not hold over a general ring $\cR$. Indeed, by moving the earring (the meridian and the arc) in \cite[Fig. 13]{kronheimer2011khovanov} to the disjoint unknot in the third picture, we obtain a skein exact triangle
\begin{equation}\label{eq: skein 1}
    \xymatrix{
\ina(Y,K,\omega;\cR)\ar[rr]^f&&\ina(Y,K,\omega;\cR)\ar[dl]\\
&\ish(Y,K,\omega;\cR)\ar[ul]&
}
\end{equation}
Hence we have for any field $\bK$,\beq\label{eq: dim K}
    \dim \ish(Y,K,\omega;\bK)=2\dim \ina(Y,K,\omega;\bK)-2\rk f\le 2\dim \ina(Y,K,\omega;\bK).
\eeq
Meanwhile, from \cite[Lemma B.1]{LY2025torsion}, we have\beq\label{eq: dim C}
    \dim \ish(Y,K,\omega;\C)\ge \dim \ina(Y,K,\omega;\C).
\eeq For more relation, see \cite{Xie2021earring}.

On the other hand, the flip symmetry in \cite[\S 2.3]{XZ2025tangle} implies that\begin{equation}\label{eq: flip symmetry}
    \ish(Y,K,\omega;\cR)\cong \ish(Y,K,\omega\cup K;\cR)
\end{equation}for the first definition because $K$ is a singular knot disjoint from the bundle set $\omega\cup \alpha$. Combining \eqref{eq: 2dim} and \eqref{eq: flip symmetry}, we obtain
\begin{equation}\label{eq: flip ina}
    \dim \ina(Y,K,\omega;\F_2)=\dim \ina(Y,K,\omega\cup K;\F_2).
\end{equation}

Finally, we study the dependence of $\omega$. Since $\omega$ is a geometric representative of Poincar\'{e} dual of the second Stiefel--Whitney class of the $SO(3)$-bundle over manifolds (or orbifolds in some approaches), the isomorphism classes of $\ina(Y,K,\omega)$ and $\ish(Y,K,\omega)$ only depend on the homology class $[\omega]\in H_1(Y\bbslash K;\F_2)$. Indeed, only the class $[\omega]\in H_1(Y;\F_2)$ matters (or $[\omega,\partial \omega]\in H_1(Y,K;\F_2)$ if $\omega$ contains arcs with endpoints on $K$); see the discussion after \cite[Proposition 3.3]{GM2023nonorientable} for the case when $Y=S^3$. 

Similar to the case of $I(Y,\omega)$ in \cite[\S 2]{LY2025dimension}, we study the concrete isomorphisms induced by (instanton) cobordism maps. We take the first approach (or the fourth approach over $\F_2$) and allow $\omega$ to contain arcs with endpoints on $K$. From \cite[\S 4.2]{kronheimer2011khovanov} (or \cite[\S 2.4]{KM2021Barnatan}), a cobordism\[(W,\Sigma,\nu):(Y_1,K_1,\omega_1)\to (Y_2,K_2,\omega_2)\]consists of 
\begin{itemize}
    \item an oriented $4$-manifold $W$ with boundary as a cobordism of $3$-manifolds $W:Y_1\to Y_2$,
    \item an possibly nonorientable surface $\Sigma$ with boundary as a cobordism of unoriented knots $\Sigma:K_1\to K_2$.
    \item a possibly nonorientable surface or a $2$-simplicial complex $\nu$ with boundary (called a \emph{bundle set} on $W$) as a cobordism of the bundle sets $\nu:\omega_1\to \omega_2$, which can have boundaries on $\Sigma$ and some transverse intersection points with $\operatorname{int}\Sigma$.
\end{itemize}

We have the following lemmas about the dependence of the bundle sets, similar to the discussion in \cite[\S 2]{LY2025dimension}.
\begin{lemma}\label{lem: cob map depending on H_2}
	Let\[(W,\Sigma,\nu_i):(Y_1,K_1,\omega_1)\to (Y_2,K_2,\omega_2)~\mathrm{for}~i=1,2\]be two cobordisms with different bundle sets satisfying\begin{equation}\label{eq: assumption nu}
	   \Sigma\cap \partial \nu_1=\Sigma\cap \partial \nu_2.
	\end{equation}Note that the transverse intersection points between $\operatorname{int}\Sigma$ and $\nu_i$ could be different. If \begin{equation}\label{eq: homology condition}
	    [\nu_1 \cup \nu_2] = 0\in H_2(W;\F_2),
	\end{equation}then for any coefficient ring $\cR$, we have
	\[
		\ina(W,\Sigma,\nu_1) = \pm \ina(W,\Sigma,\nu_2): \ina(Y_1,K_1,\omega_1; \cR) \to \ina(Y_2,K_2,\omega_2;\cR),
	\]	\[
		\aand \ish(W,\Sigma,\nu_1) = \pm \ina(W,\Sigma,\nu_2): \ish(Y_1,K_1,\omega_1; \cR) \to \ish(Y_2,K_2,\omega_2;\cR).
	\]
\end{lemma}
\brem
When $\Sigma\cap \partial\nu_1$ and $\Sigma\cap \partial \nu_2$ are not identical but have the same homology class in $H_1(W,\Sigma;\F_2)$, one might have the same results. The author does not state this stronger result because he is not sure if some similar issues as in \cite[Remark 2.2]{LY2025dimension} would happen. We will mainly focus on the cobordisms with $\Sigma\cap \partial\nu_i=\emptyset$ as $\omega_i$ do not contain arcs. Hence we do not need the strongest result.
\erem
\bpf[Proof of Lemma \ref{lem: cob map depending on H_2}]
Following \cite[\S 2.2-2.3]{kronheimer2011khovanov}, let $W^h_{\Delta}$ be obtained from $W$ by gluing $W\bbslash \Sigma$ and the double cover $\wti{\nu}_{\Delta}$ of the tubular neighborhodd $\nu=N(\Sigma)$ along $\partial \nu$ using the $2$-to-$1$ map $\partial \wti{\nu}_{\Delta}\to \partial \nu$. Here $\Delta$ is some local system on $\Sigma$ with structure group $\{\pm 1\}$, or equivalently some double cover $\Sigma_\Delta$ of $\Sigma$. There is also a non-Hausdorff space $W_\Delta$ obtained from $W\backslash \Sigma$ and $\wti{\nu}_{\Delta}$ by identifying  each point in $\wti{\nu}_{\Delta}\backslash\Sigma_\Delta$ with its image in $X\backslash \Sigma$.
 
From \cite[\S 4.2]{kronheimer2011khovanov}, we construct double cover $\Sigma_{\Delta_i}$ of $\Sigma$ for $i=1,2$ by taking $(\Sigma\backslash \omega)\times \{\pm 1\}$ and identifying across the cut with an interchange of the two sheets. In particular, if $\partial \nu_i\cap \Sigma=\emptyset$, then $\Sigma_{\Delta_i}$ is just the trivial double cover. The assumption \eqref{eq: assumption nu} implies that $\Sigma_{\Delta_1}=\Sigma_{\Delta_2}$ and hence \[W^h_{\Delta_1}=W^h_{\Delta_2}\aand W_{\Delta_1}=W_{\Delta_2}.\] Hence we omit the subscript of $\Delta_i$.

Furthermore, let $\Sigma_{\pm,i}$ be the closure of $(\Sigma\backslash\partial \nu_i)\times \{\pm 1\}$ in $W_\Delta$ and take\[\nu_i'=\pi^{-1}(\nu_i)\cup \Sigma_{-,i},\]where the map $\pi: W_\Delta\to W$. Let $\nu_i^h$ be the inverse images of $\nu_i'$ in $W_\Delta^h$. Note that a transverse intersection point of $\operatorname{int}\Sigma$ and $\nu_i$ contributes two points in  $\pi^{-1}(\nu_i)$, in which one point cancels with a point in $\Sigma_{-,i}$ (in mod $2$ sense). Hence the intersection contributes to one point in $\Sigma_{+,i}$. 

By \cite[Proposition 2.6]{kronheimer2011khovanov}, we use $\nu_i^h$ to determine singular bundle data $\bfP_i$ on $X_\Delta^h$ up to isomorphism and the addition of instantons and monopoles, which is used to construct the instanton cobordism maps. 

The assumption \eqref{eq: assumption nu} implies that $\Sigma_{-,1}=\Sigma_{-,2}$. The assumption \eqref{eq: homology condition} and \cite[Lemma 2.4 and Proposition 2.6]{kronheimer2011khovanov} imply that\[[\nu_1^h\cup \nu_2^h]=0\in H_2(W_\Delta^h;\F_2),\]and we know that the corresponding singular bundle data $\bfP_i$ are isomorphic up to the addition of instantons and monopoles. Hence the corresponding instanton cobordism maps equal up to sign.
\epf

Based on Lemma \ref{lem: cob map depending on H_2}, we have the similar result as \cite[Lemmas 2.3 and 2.5]{LY2025dimension}. The proof is also similar and we omit it.

\begin{lemma}\label{lem: I_S is an iso}
	Suppose $Y$ is a closed connected oriented $3$-manifold and suppose $K\subset Y$ is an unoriented knot. Suppose $\omega_1,\omega_1\subset Y\backslash K$ are two disjoint unoriented $1$-submanifolds such that $[\omega_1] = [\omega_2] \in H_1(Y;\F_2)$. Then there exists an embedded, possibly non-orientable surface $S\subset Y$ with $\partial S = \omega_1\cup \omega_2$. Pushing $S$ into $Y\times I$, we obtain cobordisms\[(Y\times I,K\times I,S):(Y,K,\omega_1)\to (Y,K,\omega_2)\]\[\aand (Y\times I,K\times I,S'):(Y,K,\omega_2)\to (Y,K,\omega_1).\]We write the corresponding instanton cobordism maps by $\bI_{S}$ and $\bI_{S}'$, respectively, for either $\ina$ or $\ish$. Then we have the following.
	\begin{itemize}
		\item For any choice of $S$, the maps $\bI_S$ and $\bI_{S}'$ are isomorphisms over any ring $\cR$.
		\item If $H_2(Y;\F_2) = 0$ (or equivalently $b_1(Y)=0$), then, up to sign, $\bI_S$ is independent of the choice of $S$. As a consequence, up to sign, $I(Y,\omega_1;\cR)$ and $I(Y,\omega_2;\cR)$ are canonically isomorphic.
	\end{itemize}
\end{lemma}

\section{Dimension formula}\label{sec: Dimension formula}
In this section, we prove Theorem \ref{thm: dual knot dim formula} by considering the integral case and the rational case separately. The strategy is similar, both relying on Bhat's triangle stated as follows.

\blem[{\cite[Theorem 1.1]{bhat2023newtriangle}}]\label{lem: bhat triangle}
Suppose $K\subset Y$ is a framed knot in a closed oriented connected $3$-manifold $Y$ and suppose $\omega\subset Y\backslash K$ is an unoriented $1$-submanifold. Then for any coefficient ring $\cR$, there is an exact triangle
\[
\xymatrix{
\ish(Y_0(K),\omega;\cR)\ar[rr]&&\ish(Y_2(K),\omega;\cR)\ar[dl]\\
&\ish(Y,K,\omega;\cR)\ar[ul]&
}
\]
\elem

We fix a coefficient field $\bK$ and let \[M=\nu^\sharp_{\bK}(K)\aand R=r_{\bK}(K)\]be the integers from Theorem \ref{thm: dimension formula}.

\bprop\label{prop: integral}
Suppose $K\subset S^3$ is a knot and $n\in\Z$. Then we have\[\dim \ish(S_n^3(K),\wti{K}_n;\bK)=\begin{cases}
    2R + 2|n-M|& \mathrm{if}~n\neq M;\\
    2R ~\mathrm{or}~ 2R+2&\mathrm{if}~n=M,
\end{cases}\]When $\bK=\F_2$, the case $2R+2$ will not happen and\[\dim \ina(S_n^3(K),\wti{K}_n;\F_2)=R+|n-M|.\]
\eprop
\bpf
We only consider $\ish$, and the result of $\ina$ is then obtained from \eqref{eq: 2dim}. For the mirror knot $\widebar{K}$ of $K$, we have $-S_{-n}^3(\widebar{K})\cong S^3_n(K)$ and \[\dim \ish (S^3_{-n}(\widebar{K}),\wti{\widebar{K}}_{-n};\bK)=\dim \ish (S^3_n(K),\wti{K}_n;\bK).\]By taking mirror knots, we can assume $n\ge M$. Suppose $\mu$ is the meridian of $K$ and $\lambda$ is the Seifert longitude of $K$ with $\mu\cdot \lambda=-1$. We take $(Y,K)=(S^3_n(K),\wti{K}_n)$ in Lemma \ref{lem: bhat triangle} with the meridian $\wti{\mu}=n\mu+\lambda$ and the framed longitude\[\wti{\lambda}=-\mu+k\wti{\mu}=(nk-1)\mu+k\lambda\]for some $k\in\Z$. Then we have\[\wti{\lambda}+2\wti{\mu}=(n(k+2)-1)\mu+(k+2)\lambda.\]Taking $k=-1,-2$, we obtain the exact triangles\[
\xymatrix{
\ish(S^3_{n+1}(K);\bK)\ar[rr]&&\ish(S^3_{n-1}(K);\bK)\ar[dl]\\
&\ish(S^3_n(K),\wti{K}_n;\bK)\ar[ul]&
}\]\[\xymatrix{
\ish(S^3_{(2n+1)/2}(K);\bK)\ar[rr]&&\ish(S^3;\bK)\ar[dl]\\
&\ish(S^3_n(K),\wti{K}_n;\bK)\ar[ul]&
}
\]Then we have the following dimension inequalities\begin{equation}\label{eq: dim inequal}
    \begin{aligned}
        \dim \ish(S^3_{n+1}(K);\bK)+\dim \ish(S^3_{n-1}(K);\bK)\ge &\dim \ish(S^3_n(K),\wti{K}_n;\bK),\\\dim \ish(S^3_n(K),\wti{K}_n;\bK)+\dim \ish(S^3;\bK)\ge& \dim \ish(S^3_{(2n+1)/2}(K);\bK).
    \end{aligned}
\end{equation}

We split the proof into the following three cases and apply Theorem \ref{thm: dimension formula} in each case.

\noindent{\bf Case 1}. $n\ge M+2$. We have\[\dim \ish(S^3_{n+1}(K);\bK)=R+n+1-M,~\dim \ish(S^3_{n-1}(K);\bK)=R+n-1-M,\]\[\aand \dim \ish(S^3_{(2n+1)/2}(K);\bK)=2R+2n+1-2M.\]Together with $\dim \ish(S^3;\bK)=1$, the dimension inequalities \eqref{eq: dim inequal} imply\[\dim \ish(S^3_n(K),\wti{K}_n;\bK)=2(R+n-M).\]

\noindent{\bf Case 2}. $n=M$. We have\[\dim \ish(S^3_{n+1}(K);\bK)=\dim \ish(S^3_{n-1}(K);\bK)=R+1,\]\[\aand \dim \ish(S^3_{(2n+1)/2}(K);\bK)=2R+1.\]the dimension inequalities \eqref{eq: dim inequal} imply\[2R+2\ge \dim \ish(S^3_n(K),\wti{K}_n;\bK)\ge 2R.\]The parity result in \eqref{eq: dim K} implies that the only possibilities of the dimension are $2R$ and $2R+2$. Then we exclude $2R+2$ when $\bK=\F_2$.

From the exact triangles in \cite[Lemmas 2.6 and 2.16]{LY2022integral1}, we know\[\begin{aligned}
    \chi(\ina(S^3_n(K),\wti{K}_n;\C))=&\pm \chi (\ina (S^3,K;\C)\pm \chi (\ina(S^3_{n-1}(K),\wti{K}_{n-1};\C))\\=&\chi(\ish(S^3_n(K);\C))\pmod 2
\end{aligned}
    \]where we use the identification of $\ina$ and $KHI$ from \S \ref{sec:instanton knot homology} and do not consider the Alexander garding on $KHI$. Since the Euler characteristics are independence of the coefficients, we have\[
    \chi (\ina(S^3_n(K),\wti{K}_n;\bK))=\chi (\ish(S^3_n(K);\bK))\pmod 2.\]From\[\dim \ish(S^3_{n}(K);\bK)\in \{R,R+2\},\]we know that \[\dim \ina(S^3_n(K),\wti{K}_n;\bK)\]has the same parity as $R$. Hence we again conclude the result by \eqref{eq: 2dim}. Note that the proof only works over $\F_2$ because we do not know if the parity of rank of the map $f$ in \eqref{eq: skein 1}.

\noindent{\bf Case 3}. $n=M+1$. Theorem \ref{thm: dimension formula} again implies that \[2R+4\ge \dim \ish(S^3_n(K),\wti{K}_n;\bK)\ge 2R+2.\]By adding the meridian $\mu$ of the knot $K$ to all manifolds, we also obtain\[2R+4\ge \dim \ish(S^3_n(K),\wti{K}_n,\mu;\bK)\ge 2R+2.\]Since one of \[\dim \ish(S^3_M(K);\bK)\aand \dim \ish(S^3_M(K),\mu;\bK)\]equals to $R$, we know that one of \[\dim \ish(S^3_n(K),\wti{K}_n;\bK)\aand \ish(S^3_n(K),\wti{K}_n,\mu;\bK)\]equals to $2R+2$. By \eqref{eq: flip symmetry} and the fact that $\wti{K}_n=\mu\subset S^3_n(K)$, we conclude that\[\dim \ish(S^3_n(K),\wti{K}_n;\bK)=2R+2.\]
\epf

Similar to the proof in \cite[\S 6]{LY2025dimension}, the proof of the rational case relies on the following lemma.
\begin{lemma}[{\cite[\S 4]{baldwin2020concordance} and \cite[Lemma 6.1]{LY2025dimension}}]\label{lem: existence of smaller surgery slopes}
	Suppose $p_0$ and $q_0$ are co-prime integers satisfying $p_0\neq 0$ and $|q_0|>1$. Suppose $r_0=p_0/q_0\in (k,k+1)$ for some integer $k$. Then there exist $r_i =p_i/q_i$ for $i=1,2$ that satisfy the following conditions.
	\begin{itemize}
    \item For $i=1,2$, $p_i$ and $q_i$ are co-prime, possibly zero integers, such that $p_i$ and $q_i$ have the same signs with $p_0$ and $q_0$, respectively, when they are not zero.
	\item $r_1,r_2 \in [k, k+1]$.
	\item $p_1+p_2 = p_0$ and $q_1+q_2 = q_0$.
    \item $(r_0,r_1,r_2)$ fits into a \emph{slope triad}, or more precisely i.e.\ \[p_0(-q_1)-(-p_1)q_0=(-p_1)(-q_2)-(-p_2)(-q_1)=(-p_2)q_0-p_0(-q_2)=1.\]Note that this implies that\[\frac{p_1}{q_1}>\frac{p_0}{q_0}>\frac{p_2}{q_2}.\]
	\end{itemize}
    We can further define relative prime integers $p_i,q_i$ and $r_i=p_i/q_i$ for $i=3,4$ by \[p_3=p_0+p_1,~q_3=q_0+q_1,~p_4=p_0+p_3,\aand q_4=q_0+q_3.\]Then $(r_0,r_3,r_1)$ and $(r_0,r_4,r_3)$ are also slope triads, drawn as in
    
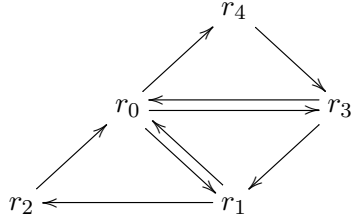
\begin{equation}\label{eq: traids of slopes}
	\xymatrix{
		&&r_4\ar[dr]&\\&r_0\ar[ur]\ar@<-2pt>[dr]\ar@<-2pt>[rr]&&r_3\ar[dl]\ar@<-2pt>[ll]\\
		r_2 \ar[ur]&& r_1\ar[ll] \ar@<-2pt>[ul]&
	}
\end{equation}
\end{lemma}
\brem
Although $n/1$ does not satisfies the assumption of Lemma \ref{lem: existence of smaller surgery slopes}, the case \[(r_0,r_1,r_2,r_3,r_4)=(\frac{n}{1},\frac{1}{0},\frac{n-1}{1},\frac{n+1}{1},\frac{2n+1}{2})\]satisfies the slope triads in \eqref{eq: traids of slopes}. Those slopes are used in the proof of Proposition \ref{prop: integral}, and motivates the choice of $r_i$ in Lemma \ref{lem: existence of smaller surgery slopes}.
\erem

\bprop\label{prop: rational}
Suppose $K\subset S^3$ is a knot and suppose $p$ and $q$ are co-prime integers with $q>1$. Then we have\[\dim \ish(S_{p/q}^3(K),\wti{K}_{p/q};\bK)=2qR + 2|p-qM|.\]\[\aand \dim \ina(S_{p/q}^3(K),\wti{K}_{p/q};\F_2)=qR + |p-qM|.\]
\eprop
\bpf
We only focus on $\ish$. The result of $\ina$ then follows from \eqref{eq: 2dim}. Again by taking mirror knots, we can assume $p/q> M$. We apply Lemma \ref{lem: existence of smaller surgery slopes} to $(p,q)=(p_0,q_0)$ to obtain the slopes $r_i$ for $i=1,2,3,4$. Then Lemma \ref{lem: bhat triangle} with suitable choices of framed knots implies the following two exact triangles.
\[
\xymatrix{
\ish(S^3_{r_3}(K);\bK)\ar[rr]&&\ish(S^3_{r_2}(K);\bK)\ar[dl]\\
&\ish(S^3_{r_0}(K),\wti{K}_{r_0};\bK)\ar[ul]&
}\]\[\xymatrix{
\ish(S^3_{r_4}(K);\bK)\ar[rr]&&\ish(S^3_{r_1}(K);\bK)\ar[dl]\\
&\ish(S^3_{r_0}(K),\wti{K}_{r_0};\bK)\ar[ul]&
}
\]Then we have dimension inequalities\begin{equation}\label{eq: dim inequal 2}
    \begin{aligned}
        \dim \ish(S^3_{r_3}(K);\bK)+\dim \ish(S^3_{r_2}(K);\bK)\ge &\ish(S^3_{r_0}(K),\wti{K}_{r_0};\bK),\\\dim\ish(S^3_{r_0}(K),\wti{K}_{r_0};\bK)+\dim \ish(S^3_{r_1}(K);\bK)\ge& \dim \ish(S^3_{r_4}(K);\bK).
    \end{aligned}
\end{equation}
From the choice of the slopes in Lemma \ref{lem: existence of smaller surgery slopes}, we have \[r_1>r_3>r_4>r_0>r_2\ge \lfloor r_0\rfloor,\]where $\lfloor r_0\rfloor$ is the maximal integer less than $r_0$. Since we assume that $r_0>M$, we know $r_2\ge M$. We split the proof into the following two cases and apply Theorem \ref{thm: dimension formula} in either case. We assume $q_i\ge 1$ for $i=1,2,3,4$.

\noindent{\bf Case 1}. $r_2>M$. We have\[\dim \ish(S^3_{r_i}(K);\bK)=q_iR+p_i-q_iM\]the dimension inequalities \eqref{eq: dim inequal 2} imply\[\dim \ish(S^3_{r_0}(K),\wti{K}_{r_0};\bK)\le (q_3R+p_3-q_3M)+(q_2R+p_2-q_2M).\]\[\dim \ish(S^3_{r_0}(K),\wti{K}_{r_0};\bK)\ge (q_4R+p_4-q_4M)-(q_1R+p_1-q_1M).\]Since\[p_4=p_0+p_3=2p_0+p_1=p_0+p_1+p_2,\]\[\aand q_4=q_0+q_3=2q_0+q_1=q_0+q_1+q_2,\]we have\[\begin{aligned}
    \dim \ish(S^3_{r_0}(K),\wti{K}_{r_0};\bK)=&(q_3R+p_3-q_3M)+(q_2R+p_2-q_2M)\\=&2(q_0R+p_0-q_0M).
\end{aligned}\]

\noindent{\bf Case 2}. $r_2=M$. Similar to Case 1, we have\[2(q_0R+p_0-q_0M)+2\ge \dim \ish(S^3_{r_0}(K),\wti{K}_{r_0};\bK)\ge 2(q_0R+p_0-q_0M).\]We use the same argument as in Case 3 in the proof of Proposition \ref{prop: integral} to conclude the result. Note that \[[\wti{K}_{p/q}]=[\mu]\in H_1(S^3_{p/q}(K);\F_2)\]and we apply Lemma \ref{lem: I_S is an iso}.
\epf

\bpf[Proof of Theorem \ref{thm: dual knot dim formula}]
The case $\omega=0$ follows directly from Propositions \ref{prop: integral} and \ref{prop: rational}. The case $\omega=\wti{K}_{p/q}$ follows from \eqref{eq: flip ina}.

For general choice of $\omega$, by Lemma \ref{lem: I_S is an iso}, we know the dimension only depends on the homology class $[\omega]\in H_1(S^3_{p/q}(K);\F_2)$. Note that 
\begin{equation}\label{eq: homology of surgery}
    H_1(S^3_{p/q}(K);\F_2) = \begin{cases}
    0 & \mathrm{when}~p~\mathrm{odd};\\
    \F_2\langle [\mu]\rangle & \mathrm{when}~p~\mathrm{even},
\end{cases}
\end{equation}where $\mu$ is the meridian of $K$. Moreover, we have $[\wti{K}_{p/q}]=[\mu]\in H_1(S^3_{p/q}(K);\F_2)$. Hence the cases $\omega=0$ and $\omega=\wti{K}_{p/q}$ imply the case of general $\omega$.
\epf
\section{\texorpdfstring{$SU(2)$}{SU(2)}-representations}\label{sec: SU(2)-representations}
In this section, we describe the relation between $\ina(Y,K)$ and $SU(2)$ representation, and then prove Theorems \ref{thm: main SU(2) 1} and \ref{thm: main SU(2) 2}.

\blem[{\cite[\S 4]{BLSY21}}]\label{lem: nondegenerate}
    Suppose $Y$ is a rational homology sphere with $H_1(Y;\Z)\cong \Z/p$ for some odd prime power $p$. Suppose $K\subset Y$ is a knot such that $H_1(Y\bbslash K;\Z)\cong \Z$. If either of the following condition holds, then\[\ina(Y,K;\Z)\cong \Z^{p}.\]
    \begin{itemize}
        \item The double branched cover $\Sigma(Y,K)$ is $SU(2)$-abelian and $\det(K)=1$.
        \item $(Y,K)$ is meridian-traceless $SU(2)$-abelian and $p$ does not divides $\det(K)$.
    \end{itemize}
In such cases, for any field $\bK$, we have \[\ina(Y,K;\bK)=\dim \ish(Y;\C)=p.\]In particular, we know $Y$ is an instanton L-space and $K$ is an instanton Floer simple knot.
\elem
\bpf
The result for the first condition follows from \cite[Theorem 4.1 and Lemma 4.4]{BLSY21}. Note that we use the notation $(Y,K)$ instead of $(L,J)$, and the conditions $H_1(Y;\Z)\cong \Z/p$ and $H_1(Y\bbslash K;\Z)\cong \Z$ imply that $[K]$ generates $H_1(Y;\Z)$, for which $K$ is called \emph{primitive}. We also replace $KHI$ with $\ina$ by the discussion in \S \ref{sec:instanton knot homology}. Indeed, we use the second approach to consider the representation variety, which is defined over $\Z$ instead of just $\C$.

The result for the second condition follows from the last paragraph in the proof of \cite[Proposition 4.5]{BLSY21}, \cite[Proposition 4.9]{BLSY21}, and the arguments in the proof of \cite[Proposition 4.1]{BLSY21} about cyclotomic polynomials.

The last equations follow from the universal coefficient theorem, \cite[Proposition 5.7]{kronheimer2011khovanov}, \cite[Theorem 1.2]{LY2020}, and \cite[Corollary 1.4]{scaduto2015instantons}, i.e.\ \[\dim \ina(Y,K;\C)=\dim KHI(Y,K)\ge \dim \ish(Y;\C)\ge |\chi (\ish(Y;\C))|=|H_1(Y;\Z)|=p.\]
\epf

\bpf[Proof of Theorems \ref{thm: main SU(2) 1} and \ref{thm: main SU(2) 2}]
We apply Lemma \ref{lem: nondegenerate} to the case\[(Y,K)=(S^3_{p/q}(K'),\wti{K}_{p/q}').\] In this case, by the proof of \cite[Lemma 4.4]{BLSY21}, we have \[H_1(Y;\Z)=|p|\aand \det(K)=|\Delta_K(-1)|=|\Delta_{K'}(-1)|=\det(K').\]By the proof of \cite[Lemma 4.2]{BLSY21}, we have\[\Sigma(S^3_{p/q}(K'),\wti{K}_{p/q}')\cong Y'_{p/2q}(\wti{K}'),\]where $\what{Y}=\Sigma(S^3,K)$ and $\what{K}'$ is the preimage of $K'$ in $\what{Y}$. Moreover, since $Y\bbslash K\cong S^3\bbslash K'$, we know that $(Y,K)$ is meridian-traceless $SU(2)$-abelian if and only if $(S^3,K')$ is $p/q$-traceless $SU(2)$-abelian.

Then it suffices to show that for $p/q\in(0,6)$, we have \[\dim \ina(S^3_{p/q}(K'),\wti{K}_{p/q}';\bK)>p\]for some field $\bK$. We take $\bK=\F_2$, then Theorems \ref{thm: dimension formula}, \ref{thm: dual knot dim formula}, and the proof of \cite[Theorem 7.3]{LY2025dimension} imply that the only possibilities are instanton L-space knots of genera at most $2$, namely $U,T_{2,3},T_{2,5}$ \cite{BLSY21,farber2024fixed}.
\epf

\bibliographystyle{alpha}

\end{document}